\documentclass[12pt]{amsart}

\usepackage{amssymb}
\usepackage{amsfonts}
\usepackage{amsmath}
\usepackage{graphicx}
\usepackage{amscd}
\usepackage{enumerate}
\usepackage{enumitem}
\newtheorem{theorem}{Theorem}
\numberwithin{theorem}{section}

\newtheorem{case}{Case}
\newtheorem{claim}{Claim}

\newtheorem{corollary}[theorem]{Corollary}

\newtheorem{definition}[theorem]{Definition}

\newtheorem{lemma}[theorem]{Lemma}

\newtheorem{question}[theorem]{Question}
\newtheorem{observation}[theorem]{Observation}
\newenvironment{proof1}[1][Proof]{\noindent\textbf{#1.} }{\ \rule{0.5em}{0.5em}}

\DeclareMathOperator{\diam}{diam}
\DeclareMathOperator{\Bd}{Bd}
\DeclareMathOperator{\ord}{ord}

\DeclareMathOperator{\Int}{int}

\begin{document}

\title[\tiny The hyperspace of non-cut subcontinua of graphs and dendrites]
{The hyperspace of non-cut subcontinua of graphs and dendrites}
\author[\tiny Hern\'andez-Guti\'errez]{Rodrigo Hern\'andez-Guti\'errez}
\author[\tiny Mart\'inez-de-la-Vega]{Ver\'onica Mart\'inez-de-la-Vega}
\author[\tiny Mart\'inez-Montejano]{Jorge M. Mart\'{\i}nez-Montejano}
\author[\tiny Vega]{Jorge E. Vega}

\address[R. Hern\'andez-Guti\'errez]{Departamento de Matem\'aticas, Universidad Aut\'onoma Metropolitana campus Iztapalapa, Av. San Rafael Atlixco 186, Leyes de Reforma 1a Secci\'on, Iztapalapa, 09310, M\'exico city, M\'exico}
\email{rod@xanum.uam.mx}

\address[V. Mart\'{\i}nez-de-la-Vega]{Instituto de Matem\'{a}ticas, Universidad Nacional Aut\'{o}noma de M\'{e}xico, Circuito Exterior, Cd. Universitaria, Ciudad de M\'{e}xico, 04510, M\'{e}xico} 
\email{vmvm@matem.unam.mx}

\address[J. M. Mart\'{\i}nez-Montejano]{Departamento de Matem\'{a}ticas, Facultad de Ciencias, Universidad Nacional Aut\'{o}noma de M\'{e}xico, Circuito Exterior, Cd. Universitaria,  Ciudad de M\'{e}xico, 04510, M\'{e}xico} \email{jorgemm@ciencias.unam.mx}

\address[J. E. Vega]{Instituto de Matem\'{a}ticas, Universidad Nacional Aut\'{o}noma de M\'{e}xico, Circuito Exterior, Cd. Universitaria, Ciudad de M\'{e}xico, 04510, M\'{e}xico} 
\email{vegacevedofc@ciencias.unam.mx}

\date{\today}

\thanks{}

\begin{abstract} 
Given a continuum $X$, let $C(X)$ denote the hyperspace of all subcontinua of $X$. In this paper we study the Vietoris hyperspace $NC^{*}(X)=\{ A \in C(X):X\setminus A\text{ is connected}\}$ when $X$ is a finite graph or a dendrite; in particular, we give conditions under which $NC^{*}(X)$ is compact, connected, locally connected or totally disconnected. Also, we prove that if $X$ is a dendrite and the set of endpoints of $X$ is dense, then $NC^{*}(X)$ is homeomorphic to the Baire space of irrational numbers.
\end{abstract}

\maketitle

\section{Introduction}
\let\thefootnote\relax\footnote{This paper was partially supported by the project ``Teor\'{\i}a de Continuos, Hiperespacios y Sistemas Din\'{a}micos III" (IN106319) of PAPIIT, DGAPA, UNAM and ``Teor\'{\i}a de Continuos e Hiperespacios, dos" (AI-S-15492) of CONACyT. Also, research of the fourth-named author was supported by the scholarship number 487053 of CONACyT.}

An \textit{arc} is a space homeomorphic to $[0,1]$. A \textit{continuum} is a compact connected metric space with more than one point. A \textit{subcontinuum} of a continuum $X$ is a nonempty closed connected subset of $X$, so one-point sets in $X$ are subcontinua of $X$. Given a continuum $X$ and a positive integer $n$, we consider the following hyperspaces:
$$2^{X}=\{A\subseteq X:A\text{ is closed and nonempty}\},$$
$$C(X)=\{A\in 2^X:A\text{ is connected}\},$$  
$$NC(X)=\{A\in C(X):X\setminus A \text{ is connected}\text{ and }\Int(A)=\emptyset \},$$
$$NC^{*}(X)=\{A\in C(X):X\setminus A \text{ is connected}\},\text{ and} $$ 
$$F_n(X)=\{A\in 2^{X}:A\text{ has at most $n$ points}\}.$$

\noindent These hyperspaces are endowed with the Hausdorff metric $H$ \cite[Definition 2.1]{Hyp} or equivalently, with the Vietoris topology \cite[Definition 1.1]{Hyp}.

The hyperspace $NC(X)$, called \textit{the hyperspace of non-cut subcontinua}, was introduced in  \cite{escobedo2017hyperspaces} inspired by the work made on \cite{Pyrih}. The purpose of this article is to study $NC^{*}(X)$ when $X$ is a finite graph or a dendrite. In Section 2 we give some preliminary definitions, results and  notation. In Section 3 we give a model of $NC^*(X)$ when $X$ is a tree and prove that, in the class of finite graphs, $NC^*(X)$ is compact if only if $X$ is either the interval or the circumference. In Section 4 we prove that if $X$ is a finite graph, then $NC^{*}(X)$ is locally connected. Finally, in section 5 we prove that if $X$ is a dendrite with its set of endpoints dense, then $NC^*(X)$ is homeomorphic to the space of irrational numbers.

\section{Preliminaries}

If $Y$ is a space and $E,F\subset Y$, then $Y=E\vert F$ means that $Y=E\cup F$, $E\cap F=\emptyset$, $E\neq\emptyset\neq F$ and both $E$ and $F$ are closed in $Y$. If $X$ is a continuum, $p,q\in X$ and $A\subset X\setminus\{p,q\}$, we say that $A$ \emph{separates} between $p$ and $q$ if there are $E,F\subset X$ such that $X\setminus A=E\vert F$, $p\in E$ and $q\in F$.

Let $(X,\rho)$ be a metric space. The open ball centered at a point $p\in X$ with radius $\varepsilon>0$ will be denoted by $B_\varepsilon^\rho(p)$.

For a given integer $n\geq 3$, a \textit{simple n-od} is a space which is
homeomorphic to the cone over an $n$-point discrete space. If $Z$ is a simple
$n$-od, then the unique point of $Z$ which is of order greater than or equal
 to 3 in $Z$ is called the core of $Z$. And a simple $3$-od is also  called
a simple triod.

A \textit{finite graph} is a continuum $X$ that can be written as the union of finitely many arcs, any two of which can intersect in at most one or both of their end points.
A \textit{tree} is a finite graph that contains no simple closed curves.

A \textit{dendroid} is an arcwise connected and hereditarily unicoherent
continuum. A \textit{dendrite} is a locally connected dendroid. We will need several properties of dendrites, all of which can be found in the compendium \cite{char}. Here we make a summary of the properties.

If $X$ is a continuum and $p\in X$ the (Menger-Urysohn) \emph{order} of $p$ in $X$, denoted by $\mathrm{ord}_X(p)$, is defined to be the least cardinal $\kappa$ such that $p$ has arbitrarily small neighborhoods with boundary of cardinality $\leq\kappa$. When $X$ is a dendrite, $\mathrm{ord}_X(p)$ is equal to the number of components of $X\setminus\{p\}$; in case this number is infinite the order is equal to the ordinal $\omega$. A point in a continuum is an \emph{endpoint} if it has order $1$, an \emph{ordinary point} if it has order $2$ and a \emph{ramification point} if it has order at least $3$. The set of endpoints (ordinary points, ramification points) is denoted by $E(X)$ ($O(X)$, $R(X)$, respectively). 

It is known that if $X$ is a dendrite, then $O(X)$ is dense in $X$, $R(X)$ is countable and $E(X)$ is of type $G_\delta$. 
\begin{lemma}\label{lemma:ramification-dense-in-arc}\cite[Theorem 4.6]{charatonik1994mapping}
For a dendrite $X$, the following are equivalent:
\begin{enumerate}[label=(\alph*)]
    \item $E(X)$ is dense,
    \item $R(X)$ is dense, and
    \item if $\alpha$ is an arc in $X$, then $\alpha\cap R(X)$ is dense.
\end{enumerate}
\end{lemma}

A dendrite $X$ is \emph{uniquely arcwise connected}, that is, for every pair $p,q\in X$ there is a unique arc $pq$ whose endpoints are exactly $p$ and $q$. Moreover, every connected subset of a dendrite is arcwise connected. If $X$ is a dendrite, $A\in C(X)$ and $p\in X\setminus A$ then there is a unique point $q\in A$ such that for every $a\in A$, $pq\subset pa$; notice that in particular $pq\cap A=\{q\}$. If $X$ is a dendrite, $p,q\in X$ and $x\in pq\setminus\{p,q\}$, it is known that $x$ separates between $p$ and $q$. 

A \textit{vertex} in a finite graph is a point of order not equal to 2. We denote by $V(X)$, the set of vertices of $X$. An \textit{edge} in a finite graph is either an arc that joins a pair of vertices and that only their end points are vertices or a circumference with at most one vertex in it. 

Let $X$ be a finite graph. Given $A$ an edge in $X$, we say that $A$ is an \textit{internal edge} if its two end-points are ramification points; otherwise, we say that $A$ is a \textit{non-internal edge}. We say that $L\subset X$ is a \textit{hair} in $X$ if $L$ is a non-internal edge with one ramification point.
We say that $L\subset X$ is a \textit{loop}
in $X$ if $L$ is a circumference and there is at most one $p\in R(X)\cap L$.

Let $X$ be a continuum. A \textit{free arc} in $X$ is any subset of $X$ of the form $A\setminus \{x,y\}$ where $A$ is an arc in $X$ with end point $x$ and $y$, and $A\setminus \{x,y\}$ is open in $X$. 

For a finite graph \(X\) and $x,y\in X$, we denote by $[x,y]$ and arc in $X$ with end points $x$ and $y$; \([x,y)=[x,y]\setminus\{y\}\); \((x,y]=[x,y]\setminus\{x\}\) and \((x,y)=[x,y]\setminus\{x,y\}\).
If $X$ is a dendrite and $a,b\in X$, there exists a unique arc in $X$ with end points $a$ and $b$. This arc is denoted $ab$.
In the case that $X$ is a tree, since $X$ is a finite graph and a dendrite, we use both notations.

A nonempty topological space is \emph{zero-dimensional} if it has a base of clopen sets (sets that are closed and open). We will use the well-known fact that a separable metric space that is the union of a countable collection of closed zero-dimensional subspaces is also zero-dimensional (see \cite[1.3.1]{Engelking}).  

A \emph{Polish} space is a separable and completely metrizable space. A space will be called \emph{nowhere compact} if it does not contain nonempty open sets with compact closures. We will use the following characterization theorem of the space of irrational numbers by Alexandroff and Urysohn.

\begin{theorem}\cite{alexandroff}\label{thm:alex-ury}
For a nonempty space $X$ the following are equivalent:
\begin{enumerate}[label=(\alph*)]
    \item $X$ is Polish, zero-dimensional and nowhere compact, and
    \item $X$ is homeomorphic to $\mathbb{R}\setminus\mathbb{Q}$.
\end{enumerate}
\end{theorem}

\section{Trees and  Finite graphs}

In this section, we present some properties of the hyperspace $NC^{*}(X)$ when $X$ is either a tree or a finite graph.  

\begin{lemma}\label{Fr_1}
If $X$ is a dendrite, $A\in NC^{*}(X)$ and $A\neq X$, then $|\Bd(A)|=1$.
\end{lemma}

\begin{proof1}
If \(\vert\Bd(A)\vert\geq2\), then $X$ contains a simple closed curve; a contradiction.
\end{proof1} 

\begin{theorem}\label{char trees}
Let $X$ be a dendrite and $A\in C(X)$. Then $A\in NC^*(X)$ if and only if one of the following holds:
\begin{enumerate}
\item $A=X$,
\item $A=\{e\}$ with $e\in E(X)$, or
\item $A=X\setminus C$ where $C$ is a component of $X\setminus \{p\}$ with $p\in X\setminus E(X)$.
\end{enumerate}
\end{theorem}

\begin{proof1}
The sufficiency is obvious; therefore, we only show the necessity.
Suppose that $A\in NC^{*}(X)$, $A\neq X$ and for each $e\in E(X)$, $A\neq\{e\}$.  We have that $\vert Bd(A)\vert=1$. Let $\{p\}=Bd(A)$. Considering that, for each $e\in E(X)$, $A\neq\{e\}$, it follows that $p\in X\setminus E(X)$. Hence, we have  that $X \setminus \{p\}=\bigcup_{i \in M}C_{i}$ with $M \subseteq\mathbb{N}$, $|M|\geq 2$ and, for each $i\in M$, $C_{i}$ is a component of $X\setminus\{p\}$. Since $p\in A$, \(X\setminus A=\bigcup_{i\in  M}C_i\setminus A\). Hence, there is $j\in M$ such that \(C_j\setminus A\neq\emptyset\). If there exists \(i\in M\), \(i\neq j\), such that \(C_i\setminus A\neq\emptyset\), then \(X\setminus A=\left(\bigcup_{i\neq j}(C_i\setminus A)\right)\mid(C_j\setminus\ A)\); a contradiction. Hence, \(\bigcup_{i\neq j}C_i\subseteq A\). If \(A\cap C_j\neq\emptyset\), then there is \(q\in\Bd(A)\cap C_j\subseteq A\setminus\{p\}\); which contradicts the fact that \(\Bd(A)=\{p\}\). We conclude that $A=X\setminus C_j$ where $C_j$ is a component of $X\setminus \{p\}$.
\end{proof1}

\begin{observation}\label{NC*(I)}
If $X$ is the interval, then $NC^{*}(X)\approx X$.
\end{observation}

\begin{observation}\label{NC*(C)}
If $X$ is a simple closed curve, then $NC^{*}(X)=C(X)$
\end{observation}

\subsection{Connectedness of $NC^*(X)$ for trees.}\
In the following result, we give a model of $NC^{*}(X)$ when $X$ is a tree.

\begin{theorem}\label{NC*(T)}
Let $X$ be a tree different from the interval $[0,1]$ with $\vert R(X)\vert=m$ and $\vert E(X)\vert=n$. Then $\textstyle NC^*(X)\approx\mathcal{K}\oplus\left(\bigoplus_{i=1}^{n+2m-2}{[0,1)}\right)$, where $\mathcal{K}$ is a simple $n$-od.
\end{theorem}

\begin{figure}[h]
\begin{center}
\includegraphics[scale=.23]{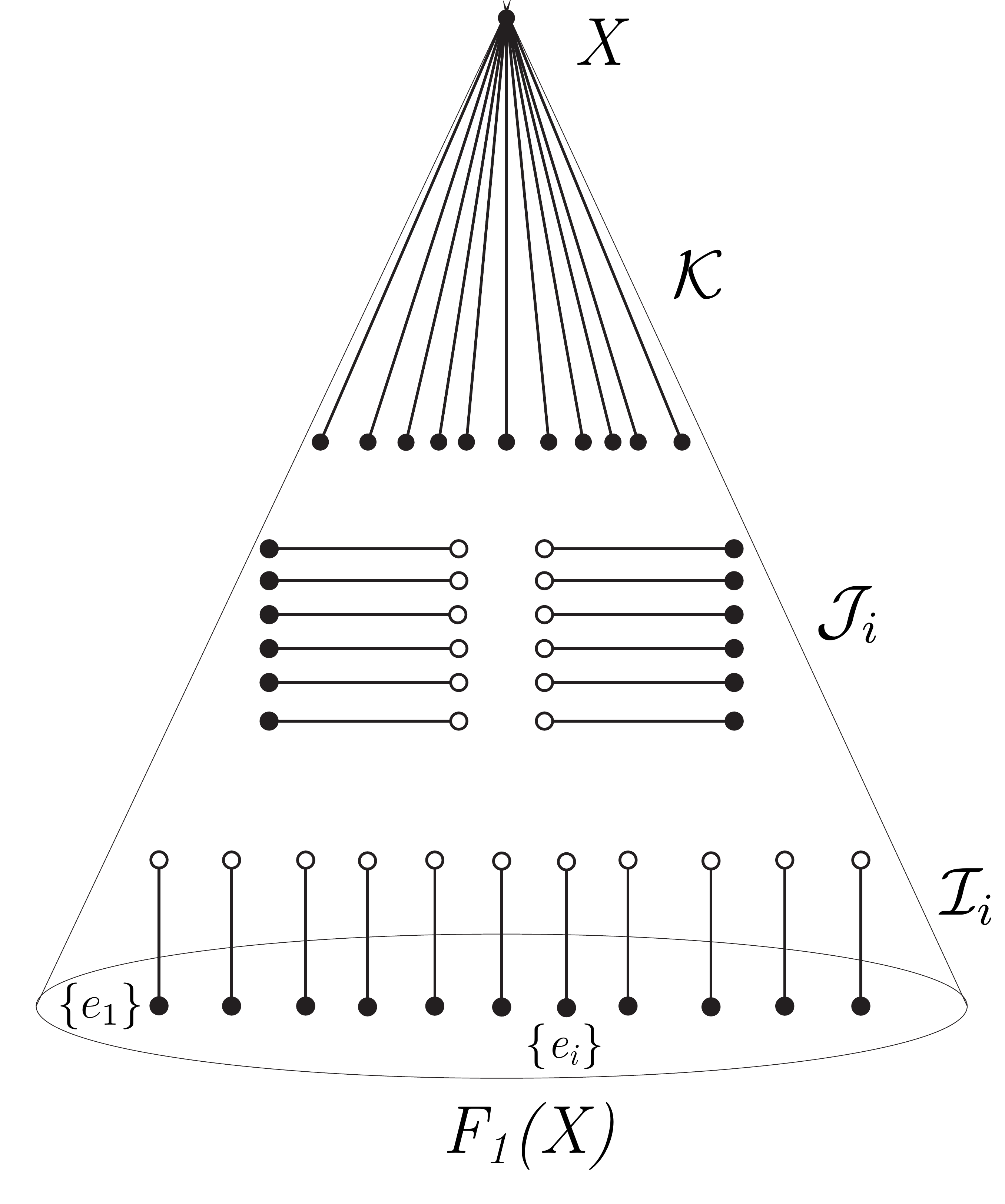}
\caption{Model of $NC^{*}(X)$ when $X$ is a tree.}
\label{NC(T)} 
\end{center}
\end{figure}

\begin{proof1}
First, notice that $X$ has $m-1$ internal edges and $n$ non-internal edges.

Let $E(X)=\{e_1,...,e_n\}$. For each $i\in\{1,\ldots,n\}$, denote by $[e_i,r_{e_i}]$  the non-internal edge that contains $e_i$. Define, for each $i\in\{1,\ldots,n\}$, $f_i:[e_i,r_{e_i})\to NC^*(X)$ by $f_i(t)=[e_i,t]$ and $g_i:[e_i,r_{e_i}]\to NC^*(X)$ by $g_i(t)=X\setminus[e_i,t)$ if $t\neq e_i$ and $g_i(e_i)=X$. It is easy to see that, for each $i\in\{1,\ldots,n\}$, both $f_i$ and $g_i$ are well defined embeddings.

Note that, given $i\in\{1,\ldots,n\}$ and $t\in(e_i,r_{e_i})$, $f_i(t)=X\setminus C_{r_{e_i}}$ where $C_{r_{e_i}}$ is  the component of $X\setminus\{t\}$
that contains $r_{e_i}$. Also, note that given $i\in\{1,\ldots,n\}$ and $t\in(e_i,r_{e_i}]$, $g_i(t)=X\setminus C_{e_i}$ where $C_{e_i}$ is the component of $X\setminus\{t\}$ that contains $e_i$.

Define $\mathcal{K}=\bigcup_{i=1}^n g_i([e_i,r_{e_i}])=\bigcup_{i=1}^{n}\{X\setminus[e_{i},t):t\in
(e_{i},r_{e_{i}}]\}\cup \{X\}$, and for each $i\in\{1,\ldots,n\}$, $\mathcal{I}_i=f_i([e_i,r_{e_i}))$. We have that $\mathcal{K}$ is a simple $n$-od (with core $X$) (see Figure \ref{NC(T)}), and for each $i\in\{1,\ldots,n\}$, $\mathcal{I}_i\approx [0,1)$ (see Figure \ref{NC(T)}) and its endpoint is $\{e_i\}$.  Given $i\in\{1,\ldots,n\}$, notice that $\mathrm{cl}_{C(X)}{(\mathcal{I}_i)}=\mathcal{I}_i\cup\{[e_i,r_{e_i}]\}$ and $[e_i,r_{e_i}]\notin NC^*(X)$ so $\mathcal{I}_i$ is closed in $NC^*(X)$.

If $m=1$, then by Theorem \ref{char trees}, $NC^\ast(X)=\mathcal{K}\cup\mathcal{I}_1\cup\cdots\cup \mathcal{I}_n$. Since $\mathcal{K},\mathcal{I}_1,\ldots, \mathcal{I}_n$ is a pairwise disjoint collection of closed sets, $NC^*(X)\approx\mathcal{K}\oplus\mathcal{I}_1\oplus\cdots\oplus\mathcal{I}_n$.

If $m\neq 1$, let us denote by $[l_1,r_1],\ldots,[l_{m-1},r_{m-1}]$ the $m-1\neq 0$ internal edges of $X$. Given $i\in\{1,\ldots,m-1\}$,  define $L_i:[l_i,r_i)\to NC^*(X)$ by $L_i(t)=X\setminus C_{r_i}$ where $C_{r_i}$ is the component of $X\setminus \{t\}$ that contains $r_i$. Also, given $i\in\{1,\ldots,m-1\}$, define $R_i:(l_i,r_i]\to NC^*(X)$ by $R_i(t)=X\setminus C_{l_i}$ where $C_{l_i}$ is the component of $X\setminus \{t\}$ that contains $l_i$. Clearly, given $i\in\{1,\ldots,m-1\}$, both $L_i$ and $R_i$ are well defined embeddings.

Define, for each $i\in\{1,\ldots,m-1\}$, $\mathcal{J}_i=L_i([l_i,r_i))$ and $\mathcal{J}_{m-1+i}=R_{i}((l_{i},r_{i}])$. We have that, for each $i\in\{1,\ldots,2m-2\}$, $\mathcal{J}_{i} \approx [0,1)$ (see Figure \ref{NC(T)}). Given $i\in\{1,\ldots,m-1\}$, $\mathrm{cl}_{C(X)}{(\mathcal{J}_i)}=\mathcal{J}_i\cup\{C\cup\{r_i\}\}$, where $C$ is the component of $X\setminus\{r_i\}$ that contains $l_i$; notice that, $C\cup\{r_i\}\notin NC^*(X)$. Likewise, given $i\in\{1,\ldots,m-1\}$, $\mathrm{cl}_{C(X)}{(\mathcal{J}_{m-1+i})}=\mathcal{J}_{m-1+i}\cup\{C\cup\{l_i\}\}$,  where $C$ is the component of $X\setminus\{l_i\}$ that contains $r_i$ and $C\cup\{l_i\}\notin NC^*(X)$. Thus, $\mathcal{J}_i$ is closed in $NC^*(X)$ for all $i\in\{1,\ldots,2m-2\}$.

Similarly to the previous case, it follows from Theorem \ref{char trees} that $NC^*(X)=\mathcal{K}\cup\mathcal{I}_1\cup\cdots\cup \mathcal{I}_n\cup\mathcal{J}_1\cup\cdots\cup\mathcal{J}_{2m-2}$ and $\mathcal{K},\mathcal{I}_1,\ldots,\mathcal{I}_n, \mathcal{J}_1,\ldots,\mathcal{J}_{2m-2}$ is a pairwise disjoint collection of closed sets. Therefore, $NC^*(X)\approx\mathcal{K}\oplus\mathcal{I}_1\oplus\cdots\oplus\mathcal{I}_n\oplus
\mathcal{J}_1\oplus\cdots\oplus\mathcal{J}_{2m-2}$ (see Figure \ref{NC(T)}).
\end{proof1}

As a consequence of Theorem \ref{NC*(T)}, we have the following three results.

\begin{corollary}\label{num_components} 
Let $X$ be a tree with $\vert R(X)\vert>0$. Then the number of components of $NC^{*}(X)$ is $2|R(X)|+|E(X)|-1$. 
\end{corollary}

\begin{corollary}\label{NC(T)_connected}
Let $X$ be a tree. Then $NC^{*}(X)$ is connected if and only if $X$ is the interval.
\end{corollary}

\begin{corollary}\label{NC(T)_continuum}
Let $X$ be a tree. $NC^{*}(X)$ is a continuum if and only if $X$ is the interval.
\end{corollary}

\subsection{Compactness of $NC^*(X)$ for graphs}\
With respect to the compactness of $NC^{*}(X)$ when $X$ is a finite graph, we have the following result.

\begin{theorem}\label{NCompacto}
Let $X$ be a finite graph. Then $NC^{*}(X)$ is compact if and only if either $X\approx[0,1]$ or $X\approx S^{1}$.
\end{theorem}

\begin{proof1}
The sufficiency is given by Propositions \ref{NC*(I)} and \ref{NC*(C)}; therefore, we only prove the necessity. Suppose that $R(X)\neq\emptyset$. In order to obtain a contradiction, we show a series of Claims.

\setcounter{claim}{0}
\begin{claim}\label{no_loops}
$X$ has no loops.
\end{claim}

Suppose to the contrary that $L$ is a loop in $X$. Let $p$ be the only ramification point in $L$. Take a sequence $\{p_n\}_{n=1}^{\infty}$ of points in $L\setminus\{p\}$ such that $p_n\to p$. We have that, for each $n\in \mathbb{N}$, $\{p_n\}\in NC^{*}(X)$, $\{p\}\not\in NC^{*}(X)$ and $\{p_n\}\to\{p\}$. Hence, $NC^*(X)$ is not compact, a contradiction. This finishes the proof of Claim \ref{no_loops}.

\begin{claim}\label{no_hairs}
$X$ has no hairs.
\end{claim}

Suppose, there is a hair $[e,r]$ in $X$ with $e\in E(X)$ and $r\in R(X)$. Note that $X\setminus \{r\}$ is not connected. In the case that $X\setminus \{r\}$ has $3$ or more components, take a sequence $\{r_n\}_{n=1}^{\infty}$ of points in $[e,r]\setminus\{r\}$ such that $r_n\to r$. We have that for each $n\in \mathbb{N}$, $[e,r_{n}]\in NC^{*}(X)$, $[e,r]\not\in NC^{*}(X)$ and $[e,r_{n}]\to[e,r]$. Hence, $NC^*(X)$ is not compact, a contradiction. Therefore, $X\setminus \{r\}$ has $2$ components, one of which is \([e,r)\); let's call the other one $K$.

Let $\mathcal{U}=\{U_n:n\in \mathbb{N}\}$ be a local basis of $X$\ at $r$ such that, for each $n\in \mathbb{N}$, $\diam(U_n)<\frac{1}{2^n}$, $\overline{U_n}$ is a simple $m$-od ($m\geq3$) with core $r$ and $\Bd(U_n)\subseteq O(X)$. Note that \(\overline{[e,r]\cup U_n}\) is also an $m$-od.

\vspace{0.2cm}
Given $n\in \mathbb{N}$, choose $z_n\in\Bd([e,r]\cup U_n)$. We show that \(\{z_n\}\in NC^*(X)$.
\vspace{0.2cm}

Since \(\ord_X(r)\geq3\), there is a point \(q\in\Bd([e,r]\cup U_n)\setminus\{z_n\}\).
Note that \(q,z_n\in K\). Since $K$ is connected, we can take an arc $[q,z_n]$ in $K$. Take a point \(w\in X\setminus\{z_n\}\). If \(w\in\overline{([e,r]\cup U_n)}\), then there is an arc \([w,q]\) in \(\overline{([e,r]\cup U_n)}\) such that \(z_n\notin[w,q]\). If \(w\notin\overline{([e,r]\cup U_n)}\), then \(w\in K\).  Hence, there is arc \([w,q]\) in $K$. Let $p$  be the first point of the arc \([w,q]\) that also belong to \([q,z_n]\). Note that, since \(q,z_{n}\in O(X)\), \(q\neq p\neq z_n\). Consider the subarcs \([w,p]\) and \([p,q]\) of \([w,q]\) and \([q,z_n]\), respectively. Thus, the arc \([w,p]\cup[p,q]\) is contained in \(K\setminus\{z_n\}\). In either case we have an arc contained in \(X\setminus\{z_n\}\) that joins \(w\) and \(q\). Therefore, $\{z_n\}\in NC^*(X)$.

So, we can construct a sequence $\{\{z_n\}\}_{n=1}^\infty$ of elements of $NC^*(X)$ such that $\{z_n\}\to\{r\}$, which contradicts that $NC^*(X)$ is compact. This concludes the proof of Claim \ref{no_hairs}.

\begin{claim}\label{Af2}
$\vert R(X)\vert\geq2$.
\end{claim}

If $\vert R(X)\vert=1$, then $X$ has either a loop or a hair, which contradicts either Claim \ref{no_loops} or Claim \ref{no_hairs}. This finishes the proof of Claim \ref{Af2}.
\vspace{0.3cm}

Let $T_{R(X)}$ be a minimal tree that contains $R(X)$.

\begin{claim}\label{Af3}
$T_{R(X)}\notin NC^*(X)$
\end{claim}

By Claim \ref{Af2}, we have that $T_{R(X)}$ is nondegenerate.  Suppose to the contrary that $T_{R(X)}\in NC^*(X)$. Then ${X\setminus T_{R(X)}}$ is an open interval $(p,q)$ in $X$ with $(p,q)\cap R(X)=\emptyset$. If $p=q$, then $\overline{X\setminus T_{R(X)}}$ is a loop in $X$, a contradiction. So, $p\neq q$. Note that $E(T_{R(X)})\neq\emptyset$ and $E(T_{R(X)})\subset R(X)$. Also, $\{p,q\}\cap E(T_{R(X)})=\emptyset$ (otherwise, if $z\in\{p,q\}\cap E(T_{R(X)})$, then $ord_X(z)=2$, a contradiction). Take $z\in E(T_{R(X)})$, since $z\notin [p,q]$, we have that $z\in int(T_{R(X)})$. Hence, $ord_X(z)=ord_{T_{R(X)}}(z)=1$; which is a contradiction. Therefore, $T_{R(X)}\notin NC^*(X)$. This concludes the proof of Claim \ref{Af3}.

\vspace{0.2cm}

Now, by Claim \ref{Af2}, $\vert R(X)\vert\geq 2$ and $T_{R(X)}$ is a nondegenerate tree. Also, by Claim \ref{Af3}, $X\setminus T_{R(X)}=\{C_1,\dots,C_m\}$ with $m\geq 2$, where ${C_i}$ is a free arc and $\overline{C_i}$ has endpoints $p_i,q_i\in R(X)\cap T_{R(X)}$. In particular, because $E(T_{R(X)})\neq\emptyset$, there exists $C\in\{C_1,\dots,C_m\}$ such that
$\overline{C}$ has endpoints $p,q$ and $p\in E(T_{R(X)})$. Note that, since $X$ has no loops, $p\neq q$. Let $A=pq$ the arc in $T_{R(X)}$ with end points $p$ and $q$.

\begin{claim}\label{Af4.1} 
$A\notin NC^{*}(X)$.
\end{claim}

If $T_{R(X)}=A$, then by Claim \ref{Af3}, $A\notin NC^{*}(X)$. So, we may assume that $T_{R(X)}\ne A$ and let $z\in T_{R(X)}\setminus A$. Since $C$ is component of $X\setminus A$ and $z\notin C$, we obtain that  $X\setminus A$ has at least two components. Therefore, $A \notin NC^*(X)$. This ends the proof of Claim 5.
\vspace{0.2cm}

Let $A\cap R(X)=\{p=p_0,\dots,q=p_k\}$ with $p_0<\dots <p_k$ in the natural order of the arc $A=pq$. Let $\alpha:[0,1]\to A$ be a homeomorphism and  let $\{s_0=0,\dots,s_k=1\}\subseteq[0,1]$
be such that $\alpha(s_i)= p_i$ for each $i\in\{0,\dots, k\}$ and define an order arc $\beta:[0,1]\rightarrow C(X)$ as $\beta(t)=\alpha([0,t])$.

\begin{claim}\label{Af4.2} 
For each $t\in[0,s_1]$,  $\beta(t)\in NC^{*}(X).$
\end{claim}

Let $0\leq t<s_1$.\newline
First, we prove that for each $z\in X\setminus \beta(t)$ there exist a point $r_z\in (X\setminus\beta(t))\cap T_{R(X)}$ and an arc $zr_z$ in $X\setminus\beta(t)$ such that $zr_z\cap T_{R(X)}=\{r_z\}$. If $z\in (X\setminus\beta(t))\cap T_{R(X)}$, let $r_z=z$. If $z\in (X\setminus\beta(t))\setminus
T_{R(X)}$, then $z\in C_i$ for some $i\in\{1,\dots,m\}$. Since $\overline {C_{i}}$ has endpoints $p_i,q_i\in T_{R(X)}\cap R(X)$, we may assume that $q_i\neq p_0$, then defining $r_z=q_i$, the arc $zr_z\subset\overline C_i$ satisfies the condition.

Now, notice that, since $\beta(0)=p_0$ and $\beta(t)\cap R(T_{R(X)})=\emptyset$, $T_{R(X)}\setminus\beta(t)$ is arcwise connected. Therefore, for each two points $u,w\in T_{R(X)}\setminus \beta(t)$, there exists an arc $uw\subset T_{R(X)}\setminus \beta(t)$. So, for each $x,y\in X\setminus\beta(t)$, $xr_x\cup r_xr_y \cup yr_y$ is a connected set in $X\setminus\beta(t)$. Hence, for $0\leq t<s_1$ the remark is proved. Since $\beta(t)\in NC^{*}(X)$ for each $t\in[0,s_1)$ and $NC^{*}(X)$ is compact, we have that $\beta(s_1)\in NC^{*}(X)$. This finishes the proof of Claim \ref{Af4.2}

\begin{claim}\label{Af4.3} Suppose that $k\geq 2$ and take $i\in\{1,\dots, k-1\}$. If $\beta(s_i)\in NC^{*}(X)$, then for each $t\in [s_i,s_{i+1}]$, $\beta(t)\in NC^{*}(X)$. 
\end{claim}

Let $s_i<t<s_{i+1}$.\newline 
First, we prove that for each $x\in X\setminus\beta(t)$ there exists an arc $L_x$ with endpoints $x,p_{i+1}$ such that $L_x\subset X\setminus \beta(t)$. Let $x\in \alpha((t,s_{i+1}] )$,  then $x=\alpha(s)$ for some $s>t$ and the arc $\alpha([s,s_{i+1}])$ satisfies the conditions. Let $x\in X\setminus\beta(s_{i+1})$ we know that there exists an arc $L$ in $X\setminus\beta(s_i)$ joining $x$ and $p_{i+1}$ since $\alpha((s_i,s_{i+1}))\cap R(X)=\emptyset$ and $p_{i}\notin L$ the first point in $xp_{i+1}\cap\alpha ([s_i,s_{i+1}])=\{p_{i+1}\}$. Thus, the arc $L=xp_{i+1}$ satisfies the condition.
\vspace{0.2cm}

Therefore, if $x,y\in X\setminus\beta(t)$, then $\{x,y\}\subseteq L_x\cup L_y$ which is a connected subset of $X\setminus\beta(t)$ and the claim is proved for $t<s_{i+1}$. Since $\beta(t)\in NC^{*}(X)$ for each $t\in[s_i,s_{i+1})$ and $NC^{*}(X)$ is compact, we conclude that  $\beta(s_{i+1})\in NC^{*}(X)$. This concludes the proof of Claim \ref{Af4.3}.

\vspace{0.2 cm}

By Claims \ref{Af4.2} and \ref{Af4.3}, we have that $A=\beta(s_k)\in NC^{*}(X)$. This contradicts Claim \ref{Af4.1}. Hence, $R(X)=\emptyset$.
\vspace{0.2 cm}

We obtain that $X$ is a finite graph with no ramification points; therefore, $X$ is homeomorphic to an interval or a simple closed curve and the theorem is proved.
\end{proof1}
\vspace{0.2 cm}

As a consequence of Theorem \ref{NCompacto}, we obtain the following result.

\begin{corollary}\label{X=NC(X)} 
Let $X$ be a finite graph. Then we have the following.
\begin{enumerate}
\item$NC^*(X)$ is a continuum if and only if $X\approx[0,1]$ or $X\approx S^1$.
\item$NC^{*}(X)\approx X$ if and only if $X$ is an interval.
\item$NC^{*}(X)= C(X)$ if and only if $X$ is a simple closed curve.
\end{enumerate}
\end{corollary}

Furthermore, we note that Theorem 10 of \cite{wilder} can be written as follows.  
\begin{theorem}\label{C(X)loccon} Let $X$ be a locally connected continuum. Then $ NC^{*}(X)= C(X)$ if and only if $X\approx S^{1}$.
\end{theorem}

\begin{question}
Let $X$ be a locally connected continuum. If $X\approx NC^{*}(X)$, then $X$ is the interval?
\end{question}

\section{Local connectedness of \(NC^*(X)\) for  finite graphs}\

We will show that if $X$ is a finite graph, then $NC^{*}(X)$ is locally connected. For this purpose we may assume that the diameter of each of the edges is one and that the metric on $X$ is convex. Now we will give a series of lemmas.

\begin{lemma}\label{uniones}
Let $X$ be a finite graph, $A\in NC^{*}(X)$  and if $A\neq X$ take a point $p_0\in X\setminus A$. Then there exists $0<\delta<\frac{1}{10}$ such that, for each $0<\varepsilon\leq\delta$, if $B\in NC^{*}(X)\cap B^H_\varepsilon(A)$, then the following conditions are satisfied.
\begin{enumerate}
\item If $A\neq X$, $p_0\notin B$.
\item If $v$ is a vertex of $X$ and $v\notin A$, then $v\notin B$.
\item If $L$ is an edge of $X$ and $L$ is not contained in $A$, then $L$ is not contained in $B$.
\item If $A\notin F_1(X)$, $L$ is an edge of $X$, $A\cap (L\setminus R(X))\neq\emptyset$ and $E$ is a component of $A\cap(L\setminus R(X))$, then $B\cap\Int(E)\neq\emptyset$; in particular, \(A\cap B\neq\emptyset\).
\item If $A\in F_1(X)$ and $A\cap B=\emptyset$, then there is an edge $L$ of $X$ such that $A\cup B\subseteq L$.
\end{enumerate}
\end{lemma}

\begin{proof1}
Let $V(X)$ be the set of vertices of $X$. If $A=X$, define \(\varepsilon=\frac{1}{10}\). If \(A\neq X\), define 
$$
 \begin{array}{l}
  \varepsilon_1 = \frac{1}{4}\min(\{d(x,y):x,y\in\Bd(A)\text{ and }x\neq y\}\cup\{1\}),\\[0.5em]
  \varepsilon_2 = \frac{1}{4}\min(\{d(x,y):x\in\Bd(A),y\in V(X)\text{ and }x\neq y\}\cup\{1\}),\\[0.5em]
  \varepsilon_3 = \min\left\{\frac{1}{4}d(p_0,A),\frac{1}{10}\right\}
  \textrm{ and } \varepsilon=\min\{\varepsilon_1,\varepsilon_2,\varepsilon_3\}.
 \end{array}
$$
If $B\in B^H_\varepsilon(A)\cap NC^*(X)$, then (1)-(5) follow. 
\end{proof1}

\begin{lemma}\label{union}
Let $X$ be a finite graph, \(A\in NC^*(X)\) and if $A\neq X$ take a point $p_0\in X\setminus A$. If $B\in C(X)$ is such that
\begin{enumerate}
\item if $A\neq X$, then $p_0\notin B$,
\item if $v$ is a vertex of $X$ and $v\notin A$, then $v\notin B$, and
\item $A\cap B\neq\emptyset$,
\end{enumerate}
then $A\cup B\in NC^*(X)$.
\end{lemma}
\begin{proof1}
If $A=X$, then the result is immediate. So, we assume that $A\neq X$. Let $q\in X\setminus(A\cup B)$ and consider an arc $l$ in $X\setminus A$ with endpoints $p_0$ and $q$. Let \(\{v_1,\ldots,v_j\}\) be the set of vertices of $X$ in $l\setminus\{p_0,q\}$ if such exist and define $v_0=p_0$ and $v_{j+1}=q$. We may assume that \(v_0<v_1<\ldots<v_{j+1}\) in the natural order of $l$; and, for each \(i\in\{0,\ldots,j\}\), let $l_i$ be the arc in $l$ joining $v_i$ and $v_{i+1}$.
We show that, for each \(i\in\{0,\ldots,j\}\), \(l_i\cap B=\emptyset\). Suppose, to the contrary, that there is $k\in\{0,\ldots,j\}$ such that \(l_k\cap B\neq \emptyset\). By (2), \(v_k,v_{k+1}\notin B\); so, since $B$ is connected, we have that \(B\subseteq l_k\subseteq l\subseteq X\setminus A\). Hence, \(A\cap B=\emptyset\); a contradiction. So, $l\cap B=\emptyset$. Thus, \(l\subseteq X\setminus (A\cup B)\). Therefore, \(A\cup B\in NC^*(X)\). 
\end{proof1}

\begin{lemma}\label{arco1}
Let $X$ be a finite graph, \(A\in NC^*(X)\), \(\delta\)
as in Lemma \ref{uniones}, $0<\varepsilon\leq\delta$ and \(B\in NC^*(X)\cap B^H_\varepsilon(A)\) such
that \(A\cap B\neq\emptyset\). Then there is a
connected subset \(\mathcal{C}\) of \(NC^*(X)\cap B^H_\varepsilon(A)\) such that \(A,A\cup B\in\mathcal{C}\).
\end{lemma}

\begin{proof1}
If either $A=X$ or $B\subseteq A$, then define \(\mathcal{C}=\{A\}\). Otherwise, take an order arc \(\beta:[0,1]\to C(X)\) such that \(\beta(0)=A\) and \(\beta(1)=A\cup B \) and define \(\mathcal{C}=\beta([0,1])\). Clearly, \(\mathcal{C}\subseteq B_\varepsilon^H(A)\). Given \(t\in[0,1]\), note that, by our choice of \(\varepsilon\), $B$ satisfies (1), (2)  and (3) of Lemma \ref{union} (for $A$). Thus, since \(\beta(t)\subseteq A\cup B\), \(\beta(t)\) satisfies (1), (2) and (3) of Lemma \ref{union} (for $A$); so, \(\beta(t)=A\cup \beta(t)\in NC^*(X)\). Therefore, \(\mathcal{C}\) is a connected subset of \(NC^*(X)\cap B^H_\varepsilon(A)\) such that \(A,A\cup B\in\mathcal{C}\).
\end{proof1}

\begin{lemma}\label{arco2}
Let $X$ be a finite graph, \(A\in NC^*(X)\), \(\delta\)
as in Lemma \ref{uniones}, $0<\varepsilon\leq\delta$ and \(B\in NC^*(X)\cap B^H_\varepsilon(A)\) such
that \(A\cap B\neq\emptyset\). Then there is a
connected subset \(\mathcal{C}\) of \(NC^*(X)\cap B^H_\varepsilon(A)\) such
that \(B,A\cup B\in\mathcal{C}\).
\end{lemma}

\begin{proof1}
If either $B=X$ or \(A\subseteq B\), then define \(\mathcal{C}=\{B\}\). Otherwise,
we construct an order arc from $B$ to $A\cup B$ as follows. First, note that, by our choice of \(\varepsilon\), the closure of each component of $A\setminus B$ is either an arc or an $n$-od. Let \(K_1,...,K_m\) be the components of $A\setminus B$. Take \(i\in\{1,\ldots,m\}\), we consider two cases.
\setcounter{case}{0}
\begin{case}
\(\overline{K_i}\) is an arc.
\end{case}
Let us write \(\overline{K_i}=[e_i,p_i]\) with \(e_i\in\Bd(B)\). Define $h_i:[0,1]\to\overline{K_i}$ by \(h_i(t)=e_i(1-t)+p_it\) and \(\beta_{i}:[0,1]\to C(X)\) by \(\beta_i(t)=B\cup h_i([0,t])\). Notice that $\beta_i$ is an order arc in $C(X)$ from $B$ to $B\cup K_i$. This ends Case 1.
\begin{case}
\(\overline{K_i}\) is an $n_i$-od.
\end{case}
Let us write $\overline{K_i}=\bigcup_{j=1}^{n_i}[e_j^i, r_i]$ where $\{e_1^i,\ldots,e_{n_i}^i\}\subseteq\Bd(B)$ and $\{r_i\}=K_i\cap R(X)$. For each $j\in\{1,\ldots,n_i\}$, define $h_j^i:[0,1]\to[e_j^i,r_i]$
by \(h_j^i(t)=e_j^i(1-t)+r_it\). Finally, define  \(\beta_{i}:[0,1]\to C(X)\) by \(\beta(t)=B\cup\bigcup_{j=1}^{n_i}h([0,t])\). Note that $\beta_i$ is an order arc in $C(X)$ from $B$ to
$B\cup K_i$. This finishes Case 2.

Now, define \(\beta:[0,1]\to C(X)\) by \(\beta(t)=\bigcup_{i=1}^m\beta_i(t)\).
We have that $\beta$ is an order arc from $B$ to $A\cup B$.
Let \(\mathcal{C}=\beta([0,1])\) and choose a point \(q_{0}\in\beta(1)\setminus\bigcup\beta([0,1))\).
Clearly, \(\mathcal{C}\subseteq
B_\varepsilon^H(A)\). Given \(t\in[0,1)\), note that, by construction, \(B \cap V(X)=\beta(t)\cap V(X)\) and $q_0\notin B\cup\beta(t)$. Hence, we can apply Lemma \ref{union} to \(B\), $q_0$ and \(\beta(t)\) to conclude that \(\beta(t)=B\cup\beta(t)\in NC^*(X)\). Also, by our choice of \(\delta\), we can apply Lemma \ref{union} to $A$ and $B$ to conclude that \(\beta(1)=A\cup B\in\ NC^*(X)\). Therefore, \(\mathcal{C}\) is a connected subset of \(NC^*(X)\cap B^H_\varepsilon(A)\)
such that \(B,A\cup B\in\mathcal{C}\).
\end{proof1}

\begin{lemma}\label{irreducible}
Let $X$ be a finite graph, $A\in NC^*(X)$, \(\delta\)
as in Lemma \ref{uniones}, \(0<\varepsilon\leq\delta\) and \(B\in NC^*(X)\cap B^H_\varepsilon(A)\) such
that \(A\cap B=\emptyset\). Then there is an edge $L$ of $X$  and an arc $K$\ such that \(A\cup B\subseteq K\subseteq L\) and \(K\in NC^*(X)\cap B_\varepsilon^H(A)\).
\end{lemma}
\begin{proof1}
Since $A\cap B=\emptyset$, by (4) of Lemma \ref{uniones}, we have that $A\in F_1(X)$. Thus, by (5) of Lemma \ref{uniones}, there is an edge $L$ of $X$ such that $A\cup B\subseteq L$. Let $K$ be the irreducible subcontinuum of $L$ such that $A\cup B\subseteq K$. Note that $K$ is an arc and $K\in B_\varepsilon^H(A)$. By our choice of \(\varepsilon\), $K$ satisfies (1), (2) and (3) of Lemma \ref{union} (for $A$); so, $K=A\cup K\in NC^*(X)$. 
\end{proof1}

\begin{theorem}
If $X$ is a finite graph, then $NC^*(X)$ is locally connected.
\end{theorem}

\begin{proof1}
Take $A\in NC^*(X)$ and $\lambda>0$. Let $\delta>0$ be a number obtained by applying Lemma \ref{uniones} to $A$ and take $0<\varepsilon\leq\min\{\delta, \lambda\}$. Let $B\in NC^*(X)$. We show that there is a connected subset \(\mathcal{C}\) of \(NC^*(X)\cap B_\varepsilon^H(A)\) such that \(A,B\in\mathcal{C}\).

In the case that \(A\cap B\neq\emptyset\), we can use Lemmas \ref{arco1}
and \ref{arco2} to construct \(\mathcal{C}\). In the case that \(A\cap B=\emptyset\), first, we use Lemma \ref{irreducible} to construct an arc \(K\in NC^*(X)\cap B_\varepsilon^H(A)\) such that \(A\cup B\subseteq K\subseteq L\); then we can apply Lemma \ref{arco1} to $A$ and $K$ to obtain a connected subset \(\mathcal{C}_1\) of  \(NC^*(X)\cap B_\varepsilon^H(A)\) such that \(A,K\in\mathcal{C}_1\). Now, take an order arc \(\alpha\) in $C(X)$ from $B$ to $K$. Let \(\mathcal{C}_2=\alpha([0,1])\) and choose a point \(q_{0}\in\alpha(1)\setminus\bigcup\alpha([0,1))\).
Clearly, \(\mathcal{C}_2\subseteq
B_\varepsilon^H(A)\). Given \(t\in[0,1)\), note that, by our choice of $\delta$, \(B
\cap V(X)=\beta(t)\cap V(X)\) and $q_0\notin B\cup\beta(t)$. Hence, we can
apply Lemma \ref{union} to \(B\), $q_0$ and \(\beta(t)\) to conclude that
\(\beta(t)=B\cup\beta(t)\in NC^*(X)\). Thus, \(\mathcal{C}=\mathcal{C}_1\cup\mathcal{C}_2\), satifies the required properties. Therefore, \(NC^*(X)\) is locally connected. 
\end{proof1}

\section{Dendrites whose set of end-points is dense}

In this section we prove that if $X$ is a dendrite with
its set of endpoints dense, then $NC^*(X)$ is homeomorphic to the space of irrational numbers; according to Theorem \ref{thm:alex-ury} we have to prove that $NC^\ast(X)$ is Polish, zero-dimensional and nowhere compact.

\subsection{Zero dimensionality}

In this subsection, we will show that that if $X$ is a dendrite with a dense set of endpoints, then $NC^{*}(X)$ is zero-dimensional. For this purpose, we will prove a series of preliminary results.

Let $X$ be a non-degenerate dendrite such that $\overline{E(X)}=X$. Fix two points $a\in X$ and $b\in X$ such that $a\neq b$. For each $x\in X\setminus (E(X)\cup\{a\})$, we define $A_{x}$ to be the component of $X\setminus\{x\}$ that contains $a$ and let $B_{x}=X\setminus A_{x}$. Next, we define
$$\mathcal{F}=\big\{B_{x}:x\in ab\setminus\{a,b\}\big\}.$$
Since for all $x\in ab\setminus \{a,b\}$,  $x\in X\setminus (E(X)\cup \{a\})$, then by Theorem \ref{char trees}(3), $\mathcal{F}\subset NC^{*}(X)$. 

\begin{observation}\label{obs1} 
If $u,w\in ab\setminus\{a,b\}$ and $w$ separates the point $u$ from the point $b$, then $A_{u}\subset A_{w}$ and $B_{w}\subset B_{u}\setminus \{u\}$.
\end{observation}

By Lemma \ref{lemma:ramification-dense-in-arc} we know that there are points $p,q\in (ab\setminus\{a,b\})\cap R(X)$ such that $q$ separates $b$ from $p$. We define:
$$\mathcal{B}=\big\{B_{x}:x\in pq\setminus \{p\}\big\}.$$
Since $p,q\in R(X)$, $X\setminus \{p\}$ and $X\setminus \{q\}$ each have at least three components, $p$ separates the point $a$ from arc $qb$ and $q$ separates arc $ap$ from point $b$. Let $C_{p}$ be a component of $X\setminus\{p\}$ that neither contains the point $a$ nor $qb\setminus\{q\}$, and $C_{q}$ a component of $X\setminus \{q\}$ that neither contains $ap\setminus\{p\}$ nor the point $b$.

\begin{observation}\label{obs2}
If $x\in pq\setminus \{p\}$, then $B_{x}\cap( A_{p}\cup C_{p})=\emptyset$ and $C_{q}\subset B_{x}$.  
\end{observation}

\begin{lemma}\label{B_abto_NC(X)}
$\mathcal{B}$ is an open set in $NC^{*}(X)$.
\end{lemma}

\begin{proof1}
Let $\mathcal{U}=\langle B_{p}\setminus\{p\}, B_{q}\setminus\overline{C_{q}}, C_{q}\rangle$. Recall that $\overline{C_{q}}=C_{q}\cup \{q\}$. We will prove that $\mathcal{U}\cap NC^{*}(X)=\mathcal{B}$.

First, let $K\in \mathcal{U}\cap NC^{*}(X)$. Since $K \subset (B_{p}\setminus\{p\})\cup(B_{q}\setminus \overline{C_{q}})\cup C_{q}$, $C_{q}\subset B_{q}$ (Observation \ref{obs2}) and $B_{q}\subset B_{p}$ (Observation \ref{obs1}) then $K\subset B_p\setminus\{p\}$. Moreover, $a\in A_{p}\subset X\setminus B_{p}$ so we have $a\notin K$ and therefore $K\neq X$. Since $K\cap (B_{q}\setminus\overline{C_{q}})\neq\emptyset\neq K\cap C_{q}$, then $q\in K$, so $K$ is not a singleton.

Since $K\in NC^{*}(X)$, by Theorem \ref{char trees}, we have $K=X\setminus C$ where $C$ is the component of $X \setminus \{y\}$ for some $y\in X\setminus E(X)$. Since $a\notin K$ and $q\in K$, we have $y\neq a$, the point $y$ separates the points $a$ and $q$, so $y\in aq\setminus\{a\}$. Since $a\notin K$, then $K=X\setminus A_{y}=B_{y}$. Notice that since $B_y=K\subset B_p\setminus\{p\}$, then $y\in B_p\setminus\{p\}$ and $y\notin ap$. So, $K=B_{y}$ for some $y\in pq\setminus\{p\}$. Therefore, $K\in\mathcal{B}$.

The other inclusion easily follows from Observations \ref{obs1} and \ref{obs2}. Therefore, $\mathcal{U}\cap NC^{*}(X)=\mathcal{B}$ and $\mathcal{B}$ is an open set in $NC^{*}(X)$.
\end{proof1}

\begin{lemma}\label{B_cerrado_NC(X)}
$\mathcal{B}$ is a closed set in $NC^{*}(X)$.
\end{lemma}

\begin{proof1}
Let $B\in\overline{\mathcal{B}}\cap NC^{*}(X)$. Then $B=\lim_{n\to\infty} B_{x_{n}}$ where $\{B_{x_{n}}\}_{n=1}^{\infty}$ is a sequence in $\mathcal{B}$. Since for all $n\in\mathbb{N}$, $x_{n}\in pq\setminus \{p\}$, we can assume that $x_{n}\to x_{0}\in pq$. So, $x_{0}\in B\cap pq$. By Observation \ref{obs2}, for all $n\in\mathbb{N}$, $A_{p}\cap B_{x_{n}}=\emptyset$ and since $A_{p}$ is an open set in $X$, we have that $A_{p}\cap B=\emptyset$. 

We have proved that $x_{0}\in B\cap pq$ and $A_{p}\cap B=\emptyset$, so that $B\neq\{e\}$ for all $e\in E(X)$ and $B\neq X$. Since $B\in NC^{*}(X)$, then by Theorem \ref{char trees}(3), $B=X\setminus C$ where $C$ is the component of $X\setminus \{ y\}$ for some $y\in X\setminus E(X)$. Since $y\in X\setminus C$ and $a\notin B$ we conclude that $y\neq a$.

If $y\neq x_0$, then $y$ separates point $a$ from point $x_0$, so $y\in ax_0\setminus\{a,x_0\}\subset ax_0\setminus\{a\}$. Otherwise, $y=x_0\in ax_0\setminus\{a\}$. Thus, $y\in ax_0\setminus\{a\}$. We will show that in fact $y\in px_{0}\setminus\{p\}$.
 
Assume to the contrary that $y\in ap\setminus\{a\}$.  Since $C_{p}\subset B_{p}$ and by Observation \ref{obs1} we have $B_{p}\subset B_{y}$, then $C_{p}\subset B_y=B$. However, by Observation \ref{obs2}, we have that for all $n\in\mathbb{N}$, $C_{p}\cap B_{x_{n}}=\emptyset$; therefore $C_{p}\cap B=\emptyset$ and we obtain a contradiction. Therefore, $y\in px_{0}\setminus\{p\}\subset pq\setminus\{p\}$ and $B=B_{y}\in \mathcal{B}$. 

We have proved that $\overline{\mathcal{B}}\cap NC^{*}(X)\subset \mathcal{B}$. Therefore, $\mathcal{B}$  is closed in $NC^{*}(X)$.
\end{proof1}

\begin{theorem}\label{dimension_cero}
Let $X$ be a dendrite such that $\overline{E(X)}=X$. Then $NC^\ast(X)$ is zero-dimensional.
\end{theorem}

\begin{proof1}
Let $\mathcal{U}=NC^{*}(X)\setminus(\{X\}\cup F_{1}(X))$. Notice that $\mathcal{U}$ is a open set in the metric space $NC^{*}(X)$, so it is an $F_\sigma$. Also, since $F_{1}(X)\cap NC^{*}(X)=F_{1}(E(X))\approx E(X)$, then $NC^{*}(X)\setminus\mathcal{U}$ is zero-dimensional. Thus, it is sufficient to show that $NC^{*}(X)$ is zero-dimensional at each point of $\mathcal{U}$.

Let $Y\in \mathcal{U}$, then by Theorem \ref{char trees}, $Y=X\setminus C$ where $C$ is a component of $X\setminus \{p\}$ and $p\notin E(X)$. Since $C\cup\{p\}$ is a non-degenerate dendrite, $\overline{C}=C\cup\{p\}$ has at least two noncut points, so there is a point $a\in E(\overline{C})\setminus\{p\}$. Notice that $a\in E(X)\cap C$. Similarly, there is a point $b\in E(X)\cap(X\setminus C)$.

Since $p$ separates the points $a$ and $b$, then $p\in ab\setminus\{a,b\}$ and using the definitions of the sets $A_{p}$ and $B_{p}$ from the begining of this section, we have that $C=A_{p}$ and $X\setminus C=X\setminus A_{p}=B_{p}$. Since $\overline{E(X)}=X$ and by Theorem \ref{lemma:ramification-dense-in-arc}, we have that there are $q\in(ap\setminus\{a,p\})\cap R(X)$ and $r\in (pb\setminus\{p,b\})\cap R(X)$.

Next, we define a countable sequence $\{\mathcal{B}_{n}:n\in \mathbb{N} \}$ of clopen sets in $NC^{*}(X)$ with $Y\in \bigcap_{n\in\mathbb{N}} \mathcal{B}_{n}$, and $\lim \diam_{H}(\mathcal{B}_{n})=0$. This will show that $Y$ is a point of zero-dimensionality of $NC^{*}(X)$. Since $p\in X\setminus E(X)$, we have the following two cases:

\setcounter{case}{0}
\begin{case}
$p \in O(X)$. 
\end{case}

Let $\{q_{n}\}_{n=1}^{\infty}$ be a sequence in $(qp\setminus \{q,p\})\cap R(X)$ with $q_{n+1}\in q_{n}p$ and $\lim q_{n}=p$. Also, let $\{r_{n}\}_{n=1}^{\infty}$ be a sequence in $(pr\setminus \{p,r\})\cap R(X)$ with $r_{n+1}\in pr _{n}$ and $\lim r_{n}=p$. For each $n\in \mathbb{N}$, let  
$$\mathcal{B}_{n}=\big\{B_{x}:x\in q_{n}r_{n}\setminus \{q_{n}\}\big\}.$$

Since $p\in q_{n}r_{n}\setminus \{q_{n}\}$ for each $n\in\mathbb{N}$, then $Y=B_{p}\in \mathcal{B}_{n}$. Notice that $\mathcal{B}_{n+1}\subset\mathcal{B}_{n}$ for each $n\in\mathbb{N}$. By Lemmas \ref{B_abto_NC(X)} and \ref{B_cerrado_NC(X)}
we know that $\mathcal{B}_n$ is clopen in $NC^{*}(X)$ for each $n\in\mathbb{N}$.

\setcounter{claim}{0}
\begin{claim} 
Given $\varepsilon >0$ there is $m\in\mathbb{N}$ such that $\diam(\mathcal{B}_{m})\leq\varepsilon$. \end{claim}

Since $\mathrm{ord}_{X}(p)=2$, there are $w_{0}\in A_{p}$ and $w_{1}\in Y\setminus\{p\}$ such that if $W$ is the component of $X\setminus \{w_{0},w_{1}\}$ that contains $p$, then $W\subset B_{\rho}(p,\frac{\varepsilon}{2})\setminus\{a,b\}$. Notice that $w_{0}\in ap$ and $w_{1}\in pb$.

Let $m\in \mathbb{N}$ be such that $q_{m},r_{m}\in W$. We show that this $m$ satisfies the Claim. Since $q_{m}\in qp\setminus \{q,p\}$ and $w_{0}\in ap$, we have that $q_{m}\in w_{0}p$. In the same way $r_{m}\in pw_{1}$. 

Take two points $x,y\in q_{m}r_{m}\setminus\{q_{m}\}$. Without loss of generality, we may assume that $y\in xr_{m}$. Given the natural order of the arc $ab$ with $a<b$, we have that: $$a< w_0 <q_{m}< x< y < r_{m}< w_{1}$$where $x,y\in W$. By Observation \ref{obs1} we have that $B_{y}\subseteq B_{x}$. We will prove that $B_{x}\subseteq B_{y}\cup W$.

Suppose to the contrary, that $B_{x}\not\subseteq B_{y}\cup W$. Then there is $z\in B_{x}$ such that $z\notin B_{y}$ and $z\notin W$.  Hence, $w_{0}\in ax$ and $X$ is a dendrite, so we have that $w_{0}$ separates between $a$ and $x$. So $X\setminus \{w_{0}\}=E\vert F$ where $a\in E$ and $x\in F$.

Now, since $B_{x}$ is connected, $w_{0}\notin B_{x}$ and $x\in B_{x}\cap F$, we have that $B_{x}\subset F$. Hence, $z\notin E$ and $z\neq w_{0}$. Secondly, since $z\notin B_{y}$, then $z\in A_{y}$. Thus, the arc $az\subset A_{y}$. If $az\cap W\neq\emptyset$, because $a,z\notin W$ and $az$ is an arc, we have that $\lvert az \cap Bd(W)\rvert\geq 2$ but since $w_{1}\notin A_{y}$, then $az \cap Bd(W)\subseteq\{w_{0}\}$, a contradiction. So, $az\cap W=\emptyset$. Since $a,x\in A_{x}\cup\{x\}$, we have that $ax\subseteq A_{x}\cup\{x\}$. In a similar way since $x,z\in B_{x}$, then $xz\subseteq B_{x}$. Hence, by the uniqueness of arcs in dendrites, $x\in az \cap W$, a contradiction. The contradiction follows from assuming that $B_{x}\not\subseteq B_{y}\cup W$, therefore $B_{x}\subseteq B_{y}\cup W$.   

From this, we get $B_{y}\subseteq B_{x}\subseteq B_{y}\cup W$. Hence, $H(B_{x},B_{y})\leq \diam(W)<\varepsilon$. This proves Claim 1.
\vspace{5pt} 

\begin{case}
$p \in R(X)$.
\end{case}

Let $\{q_{n}\}_{n=1}^{\infty}$ be a sequence in $(qp\setminus \{q,p\})\cap R(X)$ with $q_{n+1}\in q_{n}p$ and $\lim q_{n}=p$ . For each $n\in \mathbb{N}$, let
$$\mathcal{B}_{n}=\big\{B_{x}:x\in q_{n}p\setminus \{q_{n}\}\big\}.$$
Since $p\in q_{n}p\setminus \{q_{n}\}$ for each $n\in\mathbb{N}$, then $Y=B_{p}\in \mathcal{B}_{n}$. Notice that $\mathcal{B}_{n+1}\subset\mathcal{B}_{n}$ for each $n\in\mathbb{N}$. By Lemmas \ref{B_abto_NC(X)} and \ref{B_cerrado_NC(X)}
we know that $\mathcal{B}_n$ is clopen in $NC^{*}(X)$ for each $n\in\mathbb{N}$. 
 
\begin{claim}
Given $\varepsilon >0$ exists $m\in\mathbb{N}$ such that $\diam(\mathcal{B}_{m})\leq\varepsilon$. \end{claim}

Because $\overline{A_{p}}$ is a dendrite and $A_{p}=\overline{A_{p}}\setminus\{p\}$ is connected, then $\textrm{ord}_{\overline{A}_{p}}(p)=1$. Then there is $w\in \overline{A_{p}}\setminus\{p\}$ such that if $W$ is the component of $\overline{A_{p}}\setminus\{w\}$ that contains $p$, then $W\subseteq B_{\rho}(p,\frac{\varepsilon}{2})\setminus\{a\}$. Notice that $w\in ap$.

Let $m\in \mathbb{N}$ be such that $q_{m}\in W$. We show that this $m$ satisfies the Claim. Since $q_{m}\in qp\setminus \{q,p\}$ and $w\in ap$, we have that $q_{m}\in wp$. 

Take two points $x,y\in q_{m}p\setminus\{q_{m}\}$.  Without loss of generality, we may assume that $y\in xp$. Given the natural order of the arc $ab$ with $a<b$, we have that $$a< w< q_{m}< x< y < p$$ By Observation \ref{obs1}, $B_{y}\subseteq B_{x}$. We will prove that $B_{x}\subseteq B_{y}\cup W$.

Suppose to the contrary that $B_{x}\not\subseteq B_{y}\cup W$. Then there is $z\in B_{x}$ such that $z\notin B_{y}$ and $z\notin W$.  Hence, $w\in ax$ and since $X$ is a dendrite, $w$ separates between $a$ and $x$. So, $X\setminus \{w\}=E\vert F$ where $a\in E$ and $x\in F$.

Now, since $B_{x}$ is connected, $w\notin B_{x}$ and $p\in B_{x}\cap F$, we have that $B_{x}\subset F$. Hence, $z\notin E$ and $z\neq w$. Secondly, since $z\notin B_{y}$ and $Y\subset B_{y}$, we have that $z\notin Y$. Then, $a,z\in A_{p}$, so that $az\subset A_{p}$. In the dendrite $\overline{A_{p}}$, we have an arc $az$ and an open set $W$ such that $\lvert Bd_{\overline{U_{p}}}(W)\rvert=1$ and $a,z\notin W$. Hence $az\cap W=\emptyset$. Because $a,x\in A_{x}\cup\{x\}$, then $ax\subseteq A_{x}\cup\{x\}$. In similar way since $x,z\in B_{x}$, then $xz\subseteq B_{x}$. Hence, by the uniqueness arcs in dendrites, we have that $x\in az \cap W$, a contradiction. This contradiction came from assuming that $B_{x}\not\subseteq B_{y}\cup W$, therefore $B_{x}\subseteq B_{y}\cup W$.      

From this, we have that $B_{y}\subseteq B_{x}\subseteq B_{y}\cup W$. Hence, $H(B_{x},B_{y})\leq \diam(W)<\varepsilon$. This proves Claim 2.
\vspace{5pt}

Hence, $\{\mathcal{B}_{n}:n\in \mathbb{N} \}$ is a local base of clopen sets in $NC^{*}(X)$ that contains $Y$.
\end{proof1}

\subsection{Proof that $NC^*(X)$ is homeomorphic to $\mathbb{R}\setminus\mathbb{Q}$}
We start by stating the following technical result which we will need, the proof of which we leave to the reader.

\begin{lemma}\label{uniones_intersecciones}
Let $X$ be a dendrite and let $p,q\in X$ with $p\neq q$. Assume that $\{p_n\colon n\in\mathbb{N}\}\subset (pq \setminus\{p,q\})\cap O(X)$ are such that $p_{n+1}\in p_nq$ for each $n\in\mathbb{N}$ and $q=\lim_{n\to\infty}{p_n}$. For each $n\in\mathbb{N}$ let $C_n$ and $D_n$ be the components of $X\setminus\{p_n\}$ such that $p\in C_n$ and $q\in D_n$. Then

\begin{enumerate}
    \item $D_{n+1}\subset D_n$ for each $n\in\mathbb{N}$,
    \item $C_{n}\subset C_{n+1}$ for each $n\in\mathbb{N}$,
    \item $\bigcap\{\overline{D_n}\colon n\in\mathbb{N}\}=X\setminus C$ where $C$ is the component of $X\setminus\{q\}$ that contains $p$, and
    \item $\bigcup\{C_n\colon n\in\mathbb{N}\}=C$.
\end{enumerate}
\end{lemma}

\begin{theorem}\label{main}
Let $X$ be a dendrite with dense set of end-points. Then $NC^*(X)$ is homeomorphic to the space of irrational numbers.
\end{theorem}

\begin{proof1}
According to the Alexandroff-Urysohn characterization of the set of irrational numbers, we have to prove three things about $NC^*(X)$: that it is Polish, that it is zero-dimensional and that it is nowhere compact.

In \cite[Proposition 5.1]{Krupsky} it is proved that if $X$ is a locally connected continuum, the set $\mathcal{S}(X)$ of all $A\in 2^X$ such that $X\setminus A$ is disconnected is an $F_\sigma$ subset of $2^X$. Since dendrites are locally connected and $NC^\ast(X)=C(X)\setminus \mathcal{S}(X)$, it follows that $NC^\ast(X)$ is a $G_\delta$ subset of $C(X)$. As it is well-known (\cite[Lemma 1.3.12]{Engelking}), $G_\delta$ subsets of completely metrizable spaces are completely metrizable. Since $C(X)$ is compact and metrizable, $NC^\ast(X)$ is Polish.

Theorem \ref{dimension_cero} provides the zero-dimensionality. So, we are left to prove that $NC^{*}(X)$ is nowhere compact. 

Let $Y\in NC^{*}(X)$ and $\varepsilon>0$. It is sufficient to prove that the closure in $ NC^{*}(X)$ of $B_\varepsilon^{H}(Y)$ is not compact. By Theorem \ref{char trees}, we have three cases: 

\setcounter{case}{0}
\begin{case}
 $Y=\{e\}$ with $e\in E(X)$.
\end{case}

We consider $U$ an open set of $X$ such that $e\in U$ and $\diam(U)<2\varepsilon$. Since $R(X)$ is dense in $X$, there is $t\in R(X)$ such that $t\neq e$ and $t\in U$. Since $E(X)$ is dense in $X$, there is a sequence of endpoints $\{e_{k}\}_{k=1}^{\infty}$ in $U$ that converges to $t$. By Theorem \ref{char trees} (1), we know that $\{e_{k}\}\in NC^{*}(X)$ for each $k\in\mathbb{N}$. So $\{e_{k}\}\in NC^{*}(X)\cap B_\varepsilon^{H}(\{e\})$ for each $k\in\mathbb{N}$. Since $\{e_{k}\}\to\{t\}$ and $\{t\}\notin NC^{*}(X)$, we conclude that in this case $B_\varepsilon^H(Y)$ does not have compact closure in $ NC^{*}(X)$.

\begin{case}
$Y=X\setminus K$ where $K$ is a component of $X\setminus \{q\}$ with $q\in X\setminus E(X)$.
\end{case}

Let $a\in K$. Notice that $aq\cap Y=\{q\}$ and $aq\setminus\{q\}\subset K$. By Lemma \ref{lemma:ramification-dense-in-arc}, there exists $p\in aq\cap R(X)\cap B_\rho(q,\varepsilon)$. We consider a sequence $\{p_{k}\}_{k=1}^{\infty}$ of ordinary points contained in $pq$ such that $p_{k+1}\in p_{k}q$ for all $k\in \mathbb{N}$ and $p_{k}\to q$. Since $\mathrm{ord}_{X}(p_{k})=2$ for all each $k\in \mathbb{N}$, we have that $X\setminus\{p_{k}\}=C_{k}\cup D_{k}$ where $D_{k}$ is the component that contains $Y$ and $\overline{D_{k}}=D_{k}\cup \{p_{k}\}$.

By Lemma \ref{uniones_intersecciones} (1), we have that $D_{k+1}\subset D_{k}$ for each $k\in \mathbb{N}$, so that $\overline{D_{k+1}}\subset \overline{D_{k}}$ and by \cite[Exercise 4.16]{Hyp}, we conclude that $\lim_{k\to\infty} \overline{D_{k}}=\bigcap\{\overline{D_{k}}\colon k\in\mathbb{N}\}$. By Lemma \ref{uniones_intersecciones} (3), we have that $\bigcap\{\overline{D_{n}}\colon n\in\mathbb{N}\}=Y$. Now, because $Y=\lim_{k\to\infty}\overline{D_{k}}$, there is $r\in\mathbb{N}$ such that $H(Y,\overline{D_{r}})<\varepsilon$. Since $q,p_{r}\in \overline{D_{r}}$ and $\overline{D_{r}}$ is connected, we conclude that the arc $p_{r}q\subseteq\overline{D_{r}}$. 

By Lemma \ref{lemma:ramification-dense-in-arc}, $p_{r}q\cap R(X)$ is dense in $p_{r}q$. Because $p_{r}q\setminus \{p_{r},q\}$ is an open set of $p_{r}q$, there is $y\in R(X)$ such that $y\in p_{r}q\setminus \{p_{r},q\}$. Now, in the arc $yq\subset p_{r}q$ we consider a seequence $\{y_{k}\}_{k=1}^{\infty}$ of ordinary points such that $y_{k+1}\in y_{k}y$ for each $k\in \mathbb{N}$ and $y_{k}\to y$. Since $\mathrm{ord}_{X}(y_{k})=2$ for each $k\in \mathbb{N}$, we have that $X\setminus\{y_{k}\}=Y_{k}\cup X_{k}$ where $Y_{k}$ is the component that contains $Y$, does not contain $p_{r}$ and $\overline{Y_{k}}=Y_{k}\cup \{y_{k}\}$. By Theorem \ref{char trees} (3), we have that $\overline{Y_{k}}\in NC^{*}(X)$ for each $k\in\mathbb{N}$. 

By Lemma \ref{uniones_intersecciones} (2), we have that $Y_{k}\subset Y_{k+1}$ for each $k\in \mathbb{N}$, hence $\overline{Y_{k}}\subset \overline{Y_{k+1}}$ and by \cite[Exercise 4.16]{Hyp}, we conclude that 
$$\lim_{k\to\infty} \overline{Y_{k}}=\overline{\bigcup\{\overline{Y_{k}}\colon k\in\mathbb{N}\}}.$$ 
Let $Z=\lim_{k\to\infty} \overline{Y_{k}}$.

Notice that $Y\subset \overline{Y_{k}}\subset \overline{D_{r}}$ for each $k\in\mathbb{N}$. Thus, $Y\subset Z\subset\overline{D_r}$ and since $H(Y,\overline{D_{r}})<\varepsilon$, it is easy to see that $H(Y,Z)<\varepsilon$. Also, since $y_k\in \overline{Y_k}$ for all $k\in\mathbb{N}$, we obtain that $y\in Z$.

We will prove that $Z$ is not a noncut set of $X$. Since $y\in R(X)$, $X\setminus\{y\}$ has at least three components, $H$, $J$ and $K$. Without loss of generality, we may assume that $Z\setminus \{y\}\subset H$. Then $Z\cap J=\emptyset=Z\cap K$. Since $X\setminus Z\subset X\setminus\{y\}$, it follows that $J$ and $K$ are components of $X\setminus Z$. Then $X\setminus Z$ is not connected and $Z\notin NC^{*}(X)$.

Then the sequence $\{\overline{Y_{i}}\}_{k=1}^{\infty}$ of elements in $NC^{*}(X)\cap B_\varepsilon^{H}(Y)$ converges to $Z$ and $Z\notin NC^{*}(X)$. Therefore, in this case $B_\varepsilon^H(Y)$ does not have compact closure in  $NC^{*}(X)$ . 

\begin{case}
$Y=X$.  
\end{case}

Let $e\in E(X)$. Since $R(X)$ is dense in $X$, there is $p\in R(X)$ such that $p\neq e$. We consider the arc $pe$ in $X$ and $\{p_{k}\}_{k=1}^{\infty}$ a sequence of ordinary points contained in $pe$ such that $p_{k+1}\in p_{k}e$ for each $k\in \mathbb{N}$ and $p_{k}\to e$. Since $\mathrm{ord}_{X}(p_{k})=2$ for each $k\in \mathbb{N}$, we have that $X\setminus\{p_{k}\}=C_{k}\cup D_{k}$ where $C_{k}$ is the component that contains $p$ and $\overline{C_{k}}=C_{k}\cup \{p_{k}\}$.

By Lemma \ref{uniones_intersecciones} (2), we have that $C_{k}\subset C_{k+1}$ for each $k\in \mathbb{N}$, so $\overline{C_{k}}\subset \overline{C_{k+1}}$. Then by \cite[Exercise 4.16]{Hyp},  $\lim_{k\to\infty} \overline{C_{k}}=\overline{\bigcup\{\overline{C_{k}}\colon k\in\mathbb{N}\}}$. Since $e\in E(X)$, the set $X\setminus\{e\}$ is connected and by Lemma \ref{uniones_intersecciones} (4), $\bigcup\{\overline{C_{n}}\colon n\in\mathbb{N}\}=X\setminus\{e\}$; therefore, $\lim_{k\to\infty} \overline{C_{k}}=X$.

Let $r\in\mathbb{N}$ be such that $H(X,\overline{C_r})<\epsilon$. Now, since $\overline{C_r}=X\setminus D_r$, by Theorem \ref{char trees} (3), $\overline{C_r}\in NC^*(X)$. And in fact, $\overline{C_r}$ is of the second type of sets in Theorem \ref{char trees}. Then, in Case 2 we have already proved that $NC^*(X)$ is not locally compact at $\overline{C_r}$. Since $\overline{C_r}\in B_\varepsilon ^H(X)$, it follows that $B_\varepsilon^H(X)$ cannot have compact closure. Therefore, in this case $B_\varepsilon^H(Y)$ does not have compact closure in $NC^*(X)$ either.
\end{proof1}

\subsection{Some consequences}

In this section we mention some immediate consequences of our result. First, recall the following family of universal dendrites.

\begin{definition}\cite[section 2]{char}
Let $S\subset\{3,4,\ldots,\omega\}$. By $D_s$ we denote the unique dendrite satisfying the following two conditions.
\begin{enumerate}[label=(\arabic*)]
    \item If $p\in R(X)$, then $\mathrm{ord}_X(p)\in S$.
    \item For every arc $A\subset X$ and every $m\in S$ there is a point $p\in A$ such that $\mathrm{ord}_X(p)=m$.
\end{enumerate}
\end{definition}

Notice that since the order of a point in a dendrite is a topological property, it follows that if $S,T\subset\{3,4,\ldots,\omega\}$ and $D_S$ is homeomorphic to $D_T$ then $S=T$. Recall that $\mathfrak{c}=\lvert\mathbb{R}\rvert$. By taking $\mathfrak{D}=\{D_S\colon S\subset\{3,4,\ldots,\omega\}\}$ we obtain the following.

\begin{corollary}
There is a $\mathfrak{c}$-sized collection $\mathfrak{D}$ of dendrites that are pairwise non-homeomorphic such that if $X\in\mathfrak{D}$ then $NC^\ast(X)$ is homeomorphic to $\mathbb{R}\setminus\mathbb{Q}$.
\end{corollary}

Recall that a space is totally disconnected if pairs of disctinct points can be separated by clopen subsets. We can now characterize when $NC^*(X)$ is totally disconnected for a dendrite $X$.

\begin{corollary}
For a dendrite $X$ the following are equivalent:
\begin{enumerate}[label=(\alph*)]
    \item $\overline{E(X)}=X$,
    \item $NC^\ast(X)$ is homeomorphic to $\mathbb{R}\setminus\mathbb{Q}$,
    \item $NC^\ast(X)$ is zero-dimensional, and
    \item $NC^\ast(X)$ is totally disconnected.
\end{enumerate}
\end{corollary}
\begin{proof1}
(a) $\Rightarrow$ (b) follows from Theorem \ref{main}. (b) $\Rightarrow$ (c) and (c) $\Rightarrow$ (d) are obvious.

For (d) $\Rightarrow$ (a), we prove the contrapositive implication. Assume that $\overline{E(X)}\neq X$. By Lemma \ref{lemma:ramification-dense-in-arc}, there is an arc $A\subset X$ with endpoints $p$ and $q$ that has no ramification points. Since for all $x\in pq\setminus\{p,q\}$ we have that $\mathrm{ord}_{X}(x)=2$, then $X\setminus\{x\} =U_{x}\vert V_{x}$ where $p\in U_{x}$ and $q\in V_{x}$. Define $\alpha\colon pq\setminus\{p,q\}\to NC^{*}(X)$ by $\alpha(x)=A_{x}=U_{x}\cup\{x\}$. We consider $r,s \in pq\setminus \{p,q\}$ such that $r<x<s$ in the natural order of the arc $pq$ and let $\alpha\upharpoonright_{rs}\colon rs\to NC^{*}(X)$. It is easy to see that $\alpha\upharpoonright_{rs}$ is a well defined embedding.

Define $J=\alpha(rs)$. Therefore, $J\approx [0,1]$ and $J\subseteq NC^{*}(X)$. With this we conclude that $NC^{*}(X)$ is not totally disconnected.
\end{proof1} 

Now, let us turn our attention to dendroids. If $X$ is a dendroid and $x\in X$ is an endpoint in the classical sense, then clearly $\{x\}\in NC^*(X)$. Recall that there are dendroids with their set of endpoints not a $G_\delta$ (for example, see \cite[section 8]{lelek}). Also, it is known that the set of endpoints of a dendroid is a $G_\delta$ provided that the dendroid is either smooth (\cite[Theorem 3]{nikiel-tymchatyn}) or a fan (\cite[7.5]{lelek}). Thus, we may ask the following specific questions.

\begin{question}
Is there a dendroid $X$ such that $NC^*(X)$ is not a $G_\delta$ in $C(X)$?
\end{question}

\begin{question}
Let $X$ be a dendroid that is either smooth or a fan. Is $NC^*(X)$ a $G_\delta$ in $C(X)$?
\end{question}

\section*{Acknowledgements}

We would like to thank the referee for bringing the existence of Proposition 5.1 in \cite{Krupsky} to our attention; this significantly reduced the length of this paper.


\end{document}